\UseRawInputEncoding
\documentclass{amsart}
\usepackage{amsmath}
\usepackage{amssymb}
\newtheorem{thm}{Theorem}
\newtheorem{cor}[thm]{Corollary}
\newtheorem{conj}[thm]{Conjecture}
\newtheorem{lem}[thm]{Lemma}
\newtheorem{prop}[thm]{Proposition}
\newtheorem{cons}[thm]{Construction}
\theoremstyle{definition}
\newtheorem{defn}[thm]{Definition}
\theoremstyle{remark}
\newtheorem{rem}[thm]{Remark}
\newtheorem{ex}[thm]{Example}
\newtheorem{exs}[thm]{Examples}
\long\def\Thm#1{\begin{thm} #1 \end{thm}}
\long\def\Cor#1{\begin{cor} #1 \end{cor}}

\long\def\Rem#1{\begin{rem} #1 \end{rem}}
\long\def\Ex#1{\begin{ex} #1 \end{ex}}

\long\def\Ref#1#2#3#4#5#6{
\bibitem{#1}
{\rm #2,}
\textit{#3.}
{\rm #4}
\textbf{#5}
{\rm #6.}
}
\long\def\Refb#1#2#3#4{
\bibitem{#1}
{\rm #2,}
\textit{#3.}
#4.
}
\def\Rr{{\mathbb R}}

\def\phi{\varphi}

\def\leq{\leqslant}
\def\geq{\geqslant}
\def\st{\mid}

\begin{document}

\title{A Borsuk--Ulam theorem for well separated maps}

\author{M.~C.~Crabb}
\address{%
Institute of Mathematics\\
University of Aberdeen \\
Aberdeen AB24 3UE \\
UK}
\email{m.crabb@abdn.ac.uk}
\date{October 2023}
\begin{abstract}
Suppose that $f_1,\ldots ,f_m : S(V)\to\Rr$ are $m$ ($\geq 1$) continuous functions defined on the unit sphere in a Euclidean vector space $V$ 
of dimension $m+1$ satisfying $f_i(-v)=-f_i(v)$ for all $v\in S(V)$.
The classical Borsuk-Ulam theorem asserts that the image of the
map $(f_1,\ldots ,f_m) :S(V)\to\Rr^m$ contains $0=(0,\ldots ,0)$.
Pursuing ideas in \cite{bhj, frick23}, we show that a certain
separation  property will guarantee that the image contains 
an $m$-cube.
\end{abstract}
\subjclass{
52A20, 
52A38,  
55M20, 
55M25, 
57S17 
}
\keywords{Borsuk--Ulam theorem, involution, well separated}
\maketitle
Suppose that $f_1,\ldots ,f_m : S(V)\to\Rr$ are $m$ ($\geq 1$) continuous functions defined on the unit sphere in a Euclidean vector space $V$ 
of dimension $m+1$ satisfying $f_i(-v)=-f_i(v)$ for all $v\in S(V)$.
The classical Borsuk-Ulam theorem asserts that the image of the
map $(f_1,\ldots ,f_m) :S(V)\to\Rr^m$ contains $0=(0,\ldots ,0)$.
If the $m$ maps $f_i$ satisfy a certain separation  property, 
formulated below, the image will contain the $m$-cube $[-1,1]^m$.
The ideas presented here derive from the 2008 paper \cite{bhj} of
B\'ar\'any, Hubard and J\'eronimo
and the recent preprint \cite{frick23} of Frick and Wellner.
\Thm{\label{main}
Let $V$ be a Euclidean vector space of dimension $m+1>1$.
Suppose that $f_1,\ldots , f_m :S(V)\to \Rr$ are continuous functions
such that $f_i(-v)=-f_i(v)$ for all $v\in S(V)$, $i=1,\ldots ,m$,
and satisfying the condition that the open subset
$$
\Omega =\{ v\in S(V) \st |f_i(v)|<1 \text{\ for all\ }
i=1,\ldots ,m\}
$$
of $S(V)$
can be written as the disjoint union $\Omega =\Omega_+\sqcup\Omega_-$
of two open subsets which are interchanged by the antipodal involution:
$\Omega_-=-\Omega_+$.

Then the image of the continuous map
$$
(f_1,\ldots ,f_m) : S(V) \to \Rr\times\cdots\times \Rr =\Rr^m
$$
contains the $m$-cube $[-1,1]\times \cdots \times [-1,1]=[-1,1]^m$.
}
\begin{proof}
Suppose first that $|t_i|<1$ for all $i$.
The subset
$$
A_i=\{ v\in S(V) \st | f_i(v)| \leq |t_i|\}\subseteq S(V)
$$
is closed, and
$A=\bigcap_i A_i$ is contained in $\Omega$.
So there exists a continuous function $\chi : S(V)\to [-1,1]$
such that $\chi (-v)=-\chi (v)$, and
$\chi (v)=\pm 1$ if and only if 
$v\in \Omega_{\pm}\cap A$
respectively. (Indeed, knowing that the closed subsets $A_{\pm }
=\Omega_\pm\cap A$ are non-empty, by Borsuk-Ulam, we can write down
such a function $\chi$ in terms of the distance functions 
$\rho_\pm$ from $A_\pm$ -- using the standard metric, which is invariant under the antipodal involution -- as
$\chi (v)=(\rho_-(v)-\rho_+(v))/(\rho_-(v)+\rho_+(v))$.)

By the Borsuk-Ulam theorem applied to the $m$ finctions
$f_i-t_i\chi$,
there is a point $v$ such that 
$f_i(v)-t_i\chi (v)=0$ for $i=1,\ldots ,m$.
Hence $|f_i(v)|=|t_i\chi (v)|\leq |t_i|$ and thus
$v\in A_i$ for all $i$. So $|\chi (v)|=1$, and we may choose
$v$ so that $\chi (v)=1$. Then $f_i(v)=t_i$ as required.

The result for arbitrary $t_i$ follows by compactness.
\end{proof}
\Rem{When $m=1$, the condition on $\Omega$ in Theorem \ref{main}
is equivalent to the existence of a vector $v\in S(V)$ such that
$f_1(v)>1$, and so $f_1(-v)<-1$. The result is clear from the
Intermediate Value Theorem.
}
\Ex{\label{two}
Suppose that there is a hyperplane $U\subseteq V$
such that, for each $u\in S(U)$, there is some $i$ such that
$|f_i(u)|\geq 1$.
Then the condition on $\Omega$ in Theorem \ref{main} holds.
}
\begin{proof}
Choose a vector $w\in S(U^\perp)$ in the orthogonal complement of $U$. 
Then take $\Omega_+=\{v\in\Omega \, \st \,
\langle v,w\rangle >0\}$ and $\Omega_-=
\{v\in\Omega \, \st \,\langle v,w\rangle <0\}$.
\end{proof}

We can now deduce spherical versions of  the results
(a) \cite[Corollary 5.4]{frick23} and 
(b) \cite[Corollary 1]{bhj}.
Continuous functions on the unit sphere 
$S(V)$ are integrated with respect to the density given by the Euclidean metric.
\Cor{\label{three}
Suppose that $\psi_1,\ldots ,\psi_m : S(V)\to\Rr$ are continuous 
functions with $\int_{S(V)} \psi_i =1$ satisfying one of the 
conditions:
\par\noindent {\rm (a).}
There exists a hyperplane $U$ in $V$ such that,
for each $1$-dimensional subspace $L$ of $U$, 
there is some $i$ and some $u\in S(L)$,
such that $\psi_i(x)\not=0$ implies that 
$\langle x, u\rangle >0$;

\par\noindent {\rm (b).}
There exist vectors $v_I\in S(V)$ 
for each $I\subseteq \{ 1,\ldots ,m\}$ with
$v_{I'}=-v_I$,
where $I'$ denotes the complement of $I$ in $\{ 1,\ldots ,m\}$,
such that $\psi_i(x)\not=0$ for $i\in I$
implies that $\langle x,v_I\rangle > 0$.

Let $t_1,\ldots , t_m$ be real numbers, 
$-1\leq t_i \leq 1$.
Then there is a vector $v\in S(V)$ such that 
$$
\int_{\{x\in S(V)\st \langle x,v\rangle \geq 0\}} \psi_i -
\int_{\{x\in S(V)\st \langle x,v\rangle \leq 0\}} \psi_i 
=t_i
$$
for $i=1,\ldots ,m$.
}
\begin{proof}
We apply Theorem \ref{main} to the functions $f_i$ defined by
$$
f_i(v)=
\int_{\{x\in S(V)\st \langle x,v\rangle \geq 0\}} \psi_i -
\int_{\{x\in S(V)\st \langle x,v\rangle \leq 0\}} \psi_i \, .
$$
The sufficiency of condition (a) follows at once from Example
\ref{two}.

Assume that condition (b) holds and fix an orientation for $V$.
Write $U_i$ for the set of all points $x\in S(V)$ such that
$\langle x,v_I\rangle >0$ for all $I$ such that $i\in I$.
Then $U_i$ is a (non-empty) contractible open subset of the
sphere.

Any $m$ points $x_i\in U_i$, $i=1,\ldots ,m$, are linearly independent
in $V$. For suppose that $\lambda_1x_1+\ldots +\lambda_mx_m=0$. 
Put $I=\{ i \st\lambda_i>0\}$. 
Then $\langle x_i,v_I\rangle >0$ for all $i\in I$
and $\langle x_i,v_I\rangle <0$ for all $i\in I'$.
But $\sum_{i=1}^m\lambda_i\langle x_i,v_I\rangle =0$. 
So, because $\lambda_i\langle x_i,v_I\rangle \geq 0$ for all $i$,
it follows that $I=\emptyset$ and $\lambda_i=0$ for all $i$.

Now consider a point $v\in\Omega$. For each $i$,
the open hemispheres $\{ x\in S(V)\st \langle x,v\rangle >0\}$
and $\{ x\in S(V)\st \langle x,v\rangle <0\}$ must both
contain a point where $\psi_i$ is non-zero and so intersect
$U_i$. Because $U_i$ is connected,
the hyperplane $(\Rr v)^\perp$
in $V$ must, therefore,
meet each of the $m$ sets $U_i$. Choose a point $x_i
\in S((\Rr v)^\perp )\cap U_i$ for each $i=1,\ldots ,m$. 
Then $v, x_1, \ldots ,x_m$ is a basis of $V$. 
We assign $v$ to $\Omega_+$ or $\Omega_-$ according as this basis
is positively or negatively oriented.  
Since the subspace $S((\Rr v)^\perp )\cap U_i$
is contractible, this assigment does not depend on the choice of
the points $x_i$. And, because $U_i$ is open, the assigment is 
continuous.
\end{proof}
\Ex{Here are explicit examples of the two cases (a) and (b)
in Corollary \ref{three}.
Let $e_0,\ldots ,e_m$ be an orthonormal basis of $V$.
Choose $\psi_i$, concentrated near $e_i$,
such that $\psi_i(x)=0$ if $\langle x,e_i\rangle
\leq \sqrt{(m-1)/m}$.

(a). Take $U$ to be the hyperplane orthogonal to $e_0$.
If $x\in S(V)$, $\langle x,e_i\rangle > \sqrt{(m-1)/m}$, and 
$u\in S(U)$, $\langle u,e_i\rangle \geq 1/\sqrt{m}$,
then $\langle x,u\rangle >0$.

(Write $x=ae_i+y$, $u=be_i+v$, $\langle y,e_i\rangle =0$, where
$a>\sqrt{(m-1)/m}$, $b\geq 1/\sqrt{m}$, $\langle y,e_i\rangle=0$,
$\langle v,e_i\rangle =0$, $a^2+\| y\|^2=1$ and $b^2+\| v\|^2=1$.
Then $\langle x,u\rangle =ab+\langle y,v\rangle$ and
$\langle y,v\rangle^2 \leq (1-a^2)(1-b^2)=(ab)^2 -(a^2+b^2-1)$.
But $a^2+b^2>1$.)

(b). Take $v_I =(\sum_{i\in I}e_i-\sum_{i\in I'}e_i)/\sqrt{m}$.
If $x\in S(V)$, $\langle x,e_i\rangle > \sqrt{(m-1)/m}$,
and $i\in I$, so that $\langle v_I,e_i\rangle =1/\sqrt{m}$,
then $\langle x,v_I\rangle >0$, by the same calculation.
}

\end{document}